\documentclass[12pt]{article}
\usepackage{latexsym,amsmath,amssymb}

\newif\ifdviwin

\usepackage[latin1]{inputenc}
\usepackage[english]{babel}
\usepackage{indentfirst}
\usepackage[mathscr]{eucal}
\usepackage{amssymb,amsmath,amsfonts}
\usepackage{fancybox,fancyhdr}
\usepackage{graphicx}

\newif\ifdviwin

\dviwintrue 

\def\cA{\mathcal{A}}

\def\p{\psi:\Sigma\rightarrow\mathbb{L}^4}
\def\etat{\widetilde{\eta}}

\let\8=\infty \let\0=\emptyset

\let\tilde=\widetilde

\let\landa=\lambda
\let\alfa=\alpha

\let\r=\varrho

\let\parc=\partial

\def\ep{\varepsilon}
\def\vart{\vartheta}
\def\landa{\lambda}

\def\flecha{\rightarrow}
\def\esiz{\langle}
\def\esde{\rangle}

\def\cte.{\mathop{\rm cte.}\nolimits}

\def\det{\mathop{\rm det}\nolimits}

\def\Re{\mathop{\rm Re }\nolimits}

\def\cosh{\mathop{\rm cosh }\nolimits}

\def\N{\mathbb{N}}
\def\L{\mathbb{L}}

\def\M{\mathbb{M}}

\def\Q{\mathbb{Q}}
\def\R{\mathbb{R}}

\def\C{\mathbb{C}}
\def\D{\mathbb{D}}

\def\H{\mathbb{H}}
\def\S{\mathbb{S}}

\def\SL{\mathbf{SL}(2,\mathbb{C})}
\def\He{\rm{Herm}(2)}

 \newtheorem{defi}{Definition}
 \newtheorem{teo}[defi]{Theorem}
 \newtheorem{pro}[defi]{Proposition}
 \newtheorem{cor}[defi]{Corollary}
 \newtheorem{lem}[defi]{Lemma}
 
 \newtheorem{remark}[defi]{Remark}

 \newenvironment{proof}{\rm \trivlist \item[\hskip \labelsep{\it
      Proof}:]}{\par\nopagebreak \hfill $\Box$ \endtrivlist}

\numberwithin{equation}{section} 

\begin{document}
\mbox{}\vspace{0.4cm}\mbox{} 

\begin{center}
\rule{15.2cm}{1.5pt}\vspace{0.5cm} 

{\Large \bf Marginally trapped surfaces in $\L^4$ and
\\[0.3cm] an extended Weierstrass-Bryant representation}\\ \vspace{0.5cm} {\large Juan A.
Aledo$\mbox{}^a$, José A. Gálvez$\mbox{}^b$ and Pablo Mira$\mbox{}^c$}\\ 
\vspace{0.3cm} \rule{15.2cm}{1.5pt} 
\end{center}
  \vspace{1cm}

\noindent $\mbox{}^a$ Departamento de Matemáticas, Universidad de Castilla 
la Mancha. Escuela Politécnica Superior de Albacete, E-02071 Albacete, 
Spain \\ e-mail: juanangel.aledo@uclm.es \vspace{.1cm}

\noindent $\mbox{}^b$ Departamento de Geometría y Topología, Universidad 
de Granada, E-18071 Granada, Spain. \\ e-mail: 
jagalvez@ugr.es\vspace{0.1cm} 

\noindent $\mbox{}^c$ Departamento de Matemática Aplicada y Estadística, 
Universidad Politécnica de Cartagena, E-30203 Cartagena, Murcia, Spain. \\ 
e-mail: pablo.mira@upct.es \vspace{0.2cm} 

\noindent Date: February 2, 2005  

\vspace{0.3cm} 

 \begin{abstract}
We give a conformal representation in terms of meromorphic data for a 
certain class of spacelike surfaces in $\L^4$ whose mean curvature vector 
verifies $\esiz {\bf H} ,{\bf H}\esde =0$. This representation extends 
simultaneously the Weierstrass representation for minimal surfaces in 
$\R^3$ and for maximal surfaces in $\L^3$, and the Bryant representation for mean curvature one
surfaces in the hyperbolic $3$-space and in the de Sitter $3$-space.
 \end{abstract}

\section{Introduction}

In 1987 R.L. Bryant \cite{Bry} described a conformal representation for 
the surfaces with constant mean curvature $H=r$ in the hyperbolic 
$3$-space $\H^3 (-r^2)$ of constant curvature $-r^2$, similar in spirit to the classical 
Weierstrass representation for minimal surfaces in $\R^3$. Having Bryant's 
work as a basis, the theory of CMC-$r$ surfaces in $\H^3(-r^2)$ has 
experimented in the last fifteen years a great development, largely 
influenced by the global results of minimal surface theory.

Additionally, Umehara and Yamada \cite{UY} showed that there is an 
isometric deformation process via which CMC-$r$ surfaces in $\H^3(-r^2)$ 
converge analytically to minimal surfaces in $\R^3$ as $r\to 0$. This 
suggested the possibility of unifying both conformal representations into 
a more general one in a natural way.

The first such extended representation was obtained in \cite{KTUY} for a 
class of surfaces with holomorphic right Gauss map in certain Lie groups 
containing those of the form ${\bf SL} (n,\C) / {\bf SU} (n)$. An 
alternative unified representation was derived in \cite{HMN} in terms of 
Möbius geometry of surfaces.

The present paper provides a new conformal representation generalizing 
simultaneously those of Weierstrass and Bryant. To do so, we consider the 
Minkowski spacetime $\L^4$ as the ambient space, and we view 
$\R^3$ and 
$\H^3 (-r^2)$ as hyperquadrics of $\L^4$ in the usual way. 
With this, we will describe a complex representation for a certain class 
of spacelike surfaces in $\L^4$, which we will call \emph{surfaces of 
Bryant type} in $\L^4$, that includes the minimal surfaces in $\R^3$ and 
the CMC-$r$ surfaces in $\H^3 (-r^2)$. Indeed, the surfaces of Bryant type 
in $\L^4$ that lie in $\R^3\subset \L^4$ (resp. $\H^3 (-r^2) \subset 
\L^4$) are exactly the minimal surfaces of $\R^3$ (resp. the CMC-$r$ 
surfaces of 
$\H^3 (-r^2) $). 

The main geometric property of these Bryant-type surfaces is that their 
mean curvature vector $\mathbf{H}$ verifies $\esiz 
\mathbf{H},\mathbf{H}\esde =0$, where $\esiz,\esde$ is the Lorentzian 
product of $\L^4$. The spacelike surfaces defined by this condition are 
well known in General Relativity, where they are called \emph{marginally 
trapped surfaces}, and represent useful objects in the theory of 
singularities in spacetimes (see \cite{HaEl}). From our viewpoint, the 
isotropy condition 
$\esiz \mathbf{H},\mathbf{H}\esde =0$ implies that a certain \emph{Gauss 
map} of the surface in $\L^4$ is conformal, what generalizes the well 
known fact that both minimal surfaces in $\R^3$ and CMC-$r$ surfaces in 
$\H^3(-r^2)$ have conformal Gauss maps.

There are some points of special interest in the present unified conformal 
representation. First, it does not only generalize the representation 
formulae in the theories of minimal surfaces in $\R^3$ and CMC-$r$ 
surfaces in $\H^3 (-r^2)$. It also includes the conformal representations 
of their Lorentzian counterparts, namely, the theories of maximal surfaces 
in $\L^3$ \cite{Kob} and of spacelike CMC-$r$ surfaces in the de Sitter 
$3$-space $\S_1^3 (r^2)$ \cite{AiAk}, when we view $\L^3$ and $\S_1^3(r^2)$ 
as hyperquadrics of $\L^4$ in the usual way. In addition, with the present 
complex representation the Umehara-Yamada perturbation process is to some 
extent simplified, as it is viewed in the fixed ambient space $\L^4$. 
Finally, the conformal representation can be used to construct many 
complete surfaces of Bryant type in $\L^4$ which do not belong to any of 
the previous families, but that still have physical interest as they are 
marginally trapped surfaces in $\L^4$. 

The paper is organized as follows. In Section 2 we analyze the structure 
equations of a spacelike surface $\psi:\Sigma\flecha \L^4$, and prove that 
a natural \emph{hyperbolic Gauss map} $G:\Sigma\flecha \C\cup \{\8\}$ on 
the surface is conformal if $\psi$ is a marginally trapped surface. We 
also show that if the normal bundle of $\psi$ is flat, then a certain 
\emph{Hopf differential} on the surface is holomorphic. 

Section 3 describes the basic result of the present work: a conformal 
representation for the surfaces of Bryant type in $\L^4$. Here we define a 
\emph{surface of Bryant type} in $\L^4$ as a marginally trapped surface 
with flat normal bundle that is locally isometric to a minimal surface in 
$\R^3$ or to a maximal surface in $\L^3$. 

In Section 4 we will show that the meromorphic representation we have 
obtained generalizes the Weierstrass representation of the minimal (resp. 
maximal) surfaces in $\R^3$ (resp. $\L^3$), and the Bryant representation 
of the CMC-$r$ surfaces in 
$\H^3(-r^2)$ and $\S_1^3(r^2)$. We will also indicate how the 
Umehara-Yamada deformation is described in our context, and we will 
construct new examples of complete surfaces of Bryant type in $\L^4$ that 
do not belong to any of the previous families. 

Finally, in Section 5 we will classify the complete surfaces of Bryant 
type in $\L^4$ with non-negative curvature, as well as the complete 
simply-connected surfaces of Bryant type in $\L^4$ with finite total 
curvature. The paper ends up with an appendix containing some auxiliary 
results. 

\section{Marginally trapped surfaces}
Let $\L^4$ denote the $4$-dimensional Lorentz-Minkowski space, that is, 
the real vector space $\R^4$ endowed with the Lorentzian metric 
$$\esiz,\esde=-dx_0^2+dx_1^2+dx_2^2+dx_3^2,$$ in canonical coordinates. 
We shall identify $\L^4$ with the space of $2$ by 
$2$ Hermitian matrices in the usual way, 
$$(x_0,x_1,x_2,x_3)\in\L^4\longleftrightarrow 
\left(\begin{array}{cc} x_0+x_3 & x_1+i x_2\\ x_1-i x_2& x_0-x_3\\ 
\end{array}\right)\in {\rm Herm }(2).$$ Under this identification one gets 
$\esiz m,m\esde=-\det (m)$ for all $m\in {\rm Herm}(2)$. The complex Lie 
group $\SL$ acts naturally on $\L^4$ by $\Phi \cdot m =\Phi m \Phi^*$, 
being $\Phi\in \SL$, $\Phi^*=\bar{\Phi}^{t}$, and $m\in \He$. 
Consequently, $\SL$ preserves the metric and the orientations. We shall 
view the hyperbolic $3$-space of negative curvature $-r^2$ in its 
Minkowski model, that is, $\H^3(-r^2)=\{x\in \L^4: \esiz x,x\esde=-1/r^2, 
x_0>0\}.$ The above identification makes $\H^3(-r^2)$ become $$\H^3(-r^2)= 
\left\{ \frac{1}{r} \Phi \Phi^* : \Phi \in \SL\right\}, \hspace{0.5cm} 
(r>0).$$ In the same way, the \emph{de Sitter} space 
$\S_1^3(r^2)=\{x\in\L^4 ; \esiz x,x\esde =1/r^2\}$ is regarded as $$\S_1^3(r^2)= 
\left\{ \frac{1}{r} \Phi \left(\def\arraystretch{0.8}\begin{array}{cc} 0& 
1\\1& 0\\ 
\end{array}\right)\Phi^* : \Phi \in \SL\right\}, \hspace{0.5cm} 
(r>0).$$ We shall use the notation $\H^3 = \H^3(-1)$ and $\S_1^3 = \S_1^3 
(1)$. 

Finally, the \emph{positive light cone} $\N^3=\{x\in \L^4 : \esiz 
x,x\esde=0, x_0>0 \}$ is seen as the space of positive semi-definite 
matrices in $\He$ with determinant $0$, and can be described as  
$$\N^3= \left\{ w w^* : w^t=(w_1,w_2)\in \C^2\setminus\{(0,0)\}\right\},$$ where 
$w\in \C^2\setminus\{(0,0)\}$ 
is uniquely defined up to multiplication by an unimodular complex number. 
The quotient $\N^3/\R^+$ inherits a conformal structure and it can be 
regarded as the ideal boundary $\S_{\8}^2 $ of the hyperbolic $3$-space 
$\H^3$ in $\L^4$. The map $ww^*\flecha [(w_1,w_2)]$ becomes the quotient 
map of $\N^3$ onto $\S_{\8}^2$ and identifies $\S_{\8}^2$ with 
$\C\mathbf{P}^1 \equiv \C\cup \{\8\}$. 

An immersion $\psi:\Sigma\flecha\L^4$ of a connected orientable surface 
$\Sigma$ is said to be a \emph{spacelike surface} if $\Sigma$ inherits via $\psi$ a
Riemannian metric. Thus we shall regard $\Sigma$ as a Riemann surface with 
the conformal structure induced by $\psi$. 

Let $\p$ be a spacelike surface, and choose a local conformal coordinate 
$z$ on $\Sigma$ and an oriented orthonormal frame 
$\{\eta,\etat\}$ of $T^{\perp}\Sigma$, being $\etat$ a timelike vector field with values in $\H^3$, and 
$\eta$ a spacelike one. Thus the induced metric of $\Sigma$ is written as $ds^2=\landa |dz|^2$ for 
some positive smooth function $\landa$. If we define the moving frame 
\begin{equation}\label{frame}
 \sigma=(\psi_z,\psi_{\bar{z}},\eta, \etat)^T
\end{equation}
the structure equations for the immersion are 
 \begin{equation}\label{estruct}
  \sigma_z=\mathcal{U} 
\sigma,\hspace{0.5cm} \sigma_{\bar{z}}=\mathcal{V}\sigma, 
 \end{equation}
where \begin{equation*}
 \mathcal{U}= 
\left(\begin{array}{cccc} 
 (\log \landa)_z & 0 & p & \widetilde{p}\\ 
 0 & 0& E & \tilde{E}\\ 
 -2E /\landa & -2p/\landa & 0 & A \\
 2\widetilde{E}/\landa & 2\widetilde{p}/\landa & A & 0
\end{array}\right),\hspace{0.5cm} \mathcal{V}= \left(\begin{array}{cccc} 
  0 & 0& E & \tilde{E}\\ 
  0 & (\log \landa)_{\bar{z}} & \bar{p} & \bar{\tilde{p}} \\
 -2\bar{p}/\landa &  -2E /\landa & 0 & \bar{A} \\
  2\bar{\widetilde{p}}/\landa & 2\widetilde{E}/\landa & \bar{A} & 0
\end{array}\right),\end{equation*} and 
\begin{equation}\label{eaq}
 \def\arraystretch{1.2}\begin{array}{ccc} 
 E=\esiz \psi_{z\bar{z}},\eta\esde, & \tilde{E}=-\esiz \psi_{z\bar{z}},\etat\esde, & 
 A=-\esiz \eta_{z},\etat\esde,\\ 
 p=\esiz \psi_{zz},\eta\esde, & \tilde{p}=-\esiz \psi_{zz},\etat\esde . & \\ 
\end{array}
\end{equation}
The integrability condition for this system, 
$$\mathcal{U}_{\bar{z}} -\mathcal{V}_{z} 
+\left[\mathcal{U},\mathcal{V}\right] =0,$$ turns into the following 
Gauss-Codazzi-Ricci equations: 

\begin{equation}\label{GCR}\def\arraystretch{1.6} 
\begin{array}{lrcl} \text{Gauss: }& (\log \landa) _{z\bar{z}}&=& 
\frac{2}{\landa}\big(|p|^2 -|\tilde{p}|^2 +E^2 -\tilde{E}^2 \big).\\ 
 \text{Codazzi (1): } & p_{\bar{z}} -E_z &=& A\tilde{E} -\overline{A}\tilde{p}
  -E(\log \landa )_z ,\\ & \tilde{p}_{\bar{z}} -\tilde{E}_z &= &
  AE +\overline{A}p  -\tilde{E}(\log \landa )_z ,\\
   \text{Codazzi (2): } 
   & \left(\frac{\bar{p}}{\landa}\right)_z - \left(\frac{E}{\landa}\right)_{\bar{z}}
  &=& \frac{1}{\landa}\big(\bar{A}\tilde{E} -A\overline{\tilde{p}} 
  -\overline{p}(\log \landa )_z \big) \\
  & \left(\frac{\tilde{p}}{\landa}\right)_{\bar{z}} - 
  \left(\frac{\tilde{E}}{\landa}\right)_{z} &=&\frac{1}{\landa}\big(AE -\overline{A}p 
  -\tilde{p}(\log \landa )_{\bar{z}} \big) \\
 \text{Ricci: } & A_{\bar{z}} -\overline{A}_z &=& -\frac{4i}{\landa}{\rm 
 Im} (\overline{p}\tilde{p}).
\end{array}\end{equation}
The mean curvature vector of the immersion $\p$ will be denoted by 
$\mathbf{H}:\Sigma\flecha\L^4$, and with the above notations it is given by 
 \begin{equation}\label{laac}
\mathbf{H}=\frac{2}{\landa} \big( E\eta +\tilde{E}\etat\big). 
 \end{equation}
Besides, as $\eta + \etat$ takes its values in the light cone $\N^3$, we 
may define on any spacelike surface in $\L^4$ the map $[\eta + 
\etat]:\Sigma\flecha \S_{\8}^2 \equiv \C\cup \{\8\} $. It is 
straightforward to check that this map does not depend on the chosen 
orthonormal frame $\{\eta,\etat\} $ of the oriented normal bundle. So, the 
following definition makes sense: 
 \begin{defi}
The map $G=[\eta +\etat]:\Sigma\flecha \C\cup \{\8\} $ is called the 
\emph{hyperbolic Gauss map} of the spacelike surface $\psi:\Sigma\flecha 
\L^4$. 
 \end{defi}

The present paper deals with spacelike surfaces in $\L^4$ with isotropic 
mean curvature vector, that is, surfaces satisfying $\esiz 
\mathbf{H},\mathbf{H}\esde =0$. We shall call any such surface a 
\emph{marginally trapped surface}. Observe that, with this definition, any 
spacelike surface with vanishing mean curvature in $\L^4$ is marginally 
trapped. After a change of orientation in the normal bundle if necessary 
(i.e. after a change of sign in $\eta$), the above condition is written as 
$E=\tilde{E}$. Apart from their interest in Relativity Theory, the 
geometric importance of marginally trapped surfaces comes from the 
following fact. 
 \begin{lem}\label{conformaleta}
Let $\p$ be a marginally trapped surface. Then its hyperbolic Gauss map 
$G:\Sigma\flecha\C\cup\{\8\}$ is conformal. 
 \end{lem}
This result follows simply by noting that $$\esiz (\eta +\etat)_z,(\eta 
+\etat)_z \esde= 4(E-\tilde{E})(p-\tilde{p}) =0$$ for every marginally 
trapped surface in $\L^4$. 

\begin{remark}\label{casostriv}
Let $\p$ be a marginally trapped surface in $\L^4$ that actually lies in 
some $\R^3\subset \L^4$, $\L^3\subset \L^4$, $\H^3(-r^2)\subset \L^4$ or 
$\S_1^3 (r^2) \subset \L^4$. Then, by a straightforward computation, 
the condition $\esiz {\bf H},{\bf H}\esde =0$ implies that $\psi$ has 
zero mean curvature if it lies in some $\R^3$ or some $\L^3$. In the same 
way, if $\psi$ lies in some $\H^3(-r^2)$ or some $\S_1^3 (r^2)$, its mean 
curvature in that ambient space is constant, of value $H=r$. In all these 
cases, the mean curvature vector ${\bf H}$ of $\psi$ in $\L^4$ is 
parallel. 
\end{remark}

From now on we shall work with marginally trapped surfaces with \emph{flat 
normal bundle}, that is, we shall assume that the normal curvature 
vanishes identically, $R^{\perp} \equiv 0$. This amounts to say that 
$A_{\bar{z}} -\overline{A}_z =0$. But it comes clear that this condition 
implies the local existence of a real function $\beta$ on $\Sigma$ such 
that $d\beta =A dz + \bar{A} d\bar{z}$. Thus, by considering the new 
normal frame given by 
$$\xi =\cosh (\beta) \eta -\sinh (\beta) \etat,\hspace{0.5cm} \tilde{\xi} 
= -\sinh (\beta) \eta + \cosh (\beta) \etat$$ we can assume that $A=0$ 
holds in the structure equations \eqref{estruct}. In other words, there 
exists an orthonormal frame $\xi, \tilde{\xi}$ of the normal bundle that  
is parallel. Let us also remark that this parallel orthonormal frame is 
unique up to constant hyperbolic rotations in the Lorentzian normal bundle 
of the immersion. 

We shall keep denoting by $\{\eta,\etat\} $ the new parallel orthonormal 
frame $\{\xi,\tilde{\xi}\}$ of $T^{\perp}\Sigma$. It is not difficult to 
see that with the above hypothesis the equations \eqref{estruct} and 
\eqref{GCR} can be simplified considerably. 
\begin{lem}\label{flatbundle}
Let $\p$ be a marginally trapped surface with flat normal bundle, and let 
$z$ denote a local conformal coordinate on $\Sigma$. Then there is an 
orthonormal frame $\{\eta,\etat \}$ of $T^{\perp}\Sigma$ such that the 
moving frame \eqref{frame} satisfies \eqref{estruct} for the matrices 
\begin{equation}\label{matrices2}
\mathcal{U}= \left(\begin{array}{cccc} 
 (\log \landa)_z & 0 & p & \widetilde{p}\\ 
 0 & 0& E & E\\ 
 -2E /\landa & -2p/\landa & 0 & 0 \\
 2E/\landa & 2\widetilde{p}/\landa & 0 & 0
\end{array}\right),\hspace{0.5cm} \mathcal{V}= \left(\begin{array}{cccc} 
  0 & 0& E & E\\ 
  0 & (\log \landa)_{\bar{z}} & \bar{p} & \bar{\tilde{p}} \\
 -2\bar{p}/\landa &  -2E /\landa & 0 & 0 \\
  2\bar{\widetilde{p}}/\landa & 2E/\landa & 0 & 0
\end{array}\right).\end{equation} Here $E,p,\tilde{p}$ are as in \eqref{eaq}, 
and they obey the integrability conditions 
\begin{equation}\label{GaussRicci}\def\arraystretch{2} 
\begin{array}{lrcl} \text{Gauss-Ricci: }& (\log \landa) 
_{z\bar{z}}&=& \displaystyle\frac{2\big(|p|^2 -|\tilde{p}|^2)}{\landa} = 
\frac{2(\overline{p-\tilde{p}})(p+\tilde{p}) }{\landa}\\ 
 \text{Codazzi: } & p_{\bar{z}} &=& \tilde{p}_{\bar{z}} = \displaystyle\landa 
 \left(\frac{E}{\landa}\right)_z
\end{array}\end{equation}
\end{lem}
The mean curvature vector of any surface in the conditions of the above 
Lemma is given by 
 \begin{equation}\label{mincu}
  \mathbf{H}=\frac{2E}{\landa} \big(\eta +\etat\big).
  \end{equation}
It is immediate from the Codazzi equations that $\mathbf{H}$ is parallel 
if and only if $E/\landa$ is constant, if and only if both $p,\tilde{p}$ 
are holomorphic. 

Besides, let us note this other consequence of the Codazzi equations. 
 \begin{lem}
The quadratic differential $Q=(\tilde{p}-p) {\rm d}z^2$ is holomorphic on 
every marginally trapped surface with flat normal bundle in $\L^4$. 
 \end{lem} 
From now on we will denote the holomorphic function $\tilde{p}-p$ as 
$q=\tilde{p}-p$. We shall call $Q = q(z) dz^2$ the \emph{Hopf 
differential} of the surface. 
\begin{remark}\label{fl}
If $\p$ is a marginally trapped surface with flat normal bundle on which 
$Q$ vanishes identically, then by \eqref{matrices2}, $\eta+\etat$ is constant and $\Sigma$ is 
flat. Thus $\psi(\Sigma)$ lies in an affine degenerate hyperplane of 
$\L^4$. Flat surfaces lying in degenerate hyperplanes of $\L^4$ were 
completely described in explicit coordinates in \cite{GMM1}. Let us also 
indicate that many of these surfaces are complete. 
\end{remark}
In the remaining of the present work we will assume that $Q$ does not 
vanish identically.

\section{A conformal representation}
Let us start this section recalling a classical result by Ricci, stating 
that a necessary and sufficient condition for a Riemannian surface 
$(S,ds^2)$ to be locally isometric to a minimal surface in $\R^3$ is that 
it has non-positive curvature $K\leq 0$ and the conformal pseudo-metric 
$d\tilde{s}^2 =\sqrt{-K}ds^2$ is flat in case $K\not\equiv 0$. This is 
equivalent to impose that the conformal pseudo-metric $-K ds^2$ has 
constant curvature $1$ at its regular points. 

Analogously, it can be easily proved that $(S,ds^2)$ is locally isometric 
to a maximal surface in $\L^3$ if and only if $K\geq 0$ and $\sqrt{K} 
ds^2$ is flat, if and only if $K\geq 0$ and $K ds^2$ has constant 
curvature $-1$.
\begin{defi}
Let $\p$ be a marginally trapped surface with flat normal bundle. Then 
$\psi$ is said to be a \emph{surface of Bryant type} in $\L^4$ provided 
$\Sigma$ is locally isometric to some minimal surface in $\R^3$ or to some 
maximal surface in $\L^3$. 
\end{defi}
Let us examine this kind of surfaces. For this we start with a simply 
connected marginally trapped surface with flat normal bundle $\p$. From 
the Gauss equation in \eqref{GaussRicci} we find that its Gauss curvature 
is $$K= \frac{4\left(|\tilde{p}|^2-|p|^2\right)}{\landa^2}.$$ We shall 
denote $\varepsilon={\rm sign} (K)={\rm sign} 
\left(|\tilde{p}|^2-|p|^2\right)$. If $\varepsilon=0$ the surface is of 
Bryant type. Otherwise we can define on $\Sigma$ away from the isolated 
flat points of $\Sigma$ the metric $d\tilde{s}^2=\sqrt{\varepsilon K} 
ds^2$, and the above formula easily implies that $d\tilde{s}^2$ is flat if 
and only if 
$$\left(\log \left(\varepsilon (|\tilde{p}|^2-|p|^2)\right) 
\right)_{z\bar{z}} =0,$$ that is, if and only if there is a holomorphic 
function $\varphi:\Sigma\flecha \C$ such that $\varepsilon 
(|\tilde{p}|^2-|p|^2) =|\varphi|^2$ (recall that $\Sigma$ is simply 
connected). But by the Ricci equation, this is equivalent to the fact that 
\begin{equation}\label{defalfa}
 \tilde{p}+p = \varepsilon |f|^2 q,
 \end{equation}
being 
\begin{equation}\label{masalfa}
f= \varphi/ q:\Sigma\flecha \C\cup \{\8\} 
\end{equation} a meromorphic function on $\Sigma$.

\begin{remark}\label{rem7}
As we know that a marginally trapped surface with flat normal bundle has 
parallel mean curvature if and only if both $p,\tilde{p}$ are holomorphic, 
it turns out that a surface of Bryant type has parallel mean curvature if 
and only if $f$ is constant. 

Moreover, if $\p$ is a marginally trapped surface with flat normal bundle 
and parallel mean curvature vector in $\L^4$, then it belongs to a known 
family. Indeed, if this is the case we get from the Ricci equation that 
$\widetilde{p}=k p$, $k\in \R$, which, jointly with the structure equations (\ref{matrices2}) 
provide $\psi_z=k'(k\eta+\widetilde{\eta})_z$, $k'\in\R$. This identity
implies that $\psi$ must lie in a totally umbilical affine hyperquadric of 
$\L^4$, and thus is a known example by Remark \ref{casostriv}.
\end{remark}

The basic result of the present work is a conformal representation for 
surfaces of Bryant type in $\L^4$. 

\begin{teo}[Conformal representation]\label{rep}
Let $\psi:\Sigma\flecha\L^4$ be a non-flat simply connected surface of 
Bryant type, and let $\varepsilon=\pm 1$ denote the sign of its Gaussian 
curvature. Then there exist a meromorphic function $g$ and a holomorphic 
$1$-form $\omega$ on $\Sigma$, and three constants $a,b\in \R$, $c\in \C$
satisfying 
\begin{enumerate} 
\item[{\bf C.1}] $1-\varepsilon |g|^2>0$, and the zeros of 
$\omega$ of order $2k$ correspond to the poles of $g$ of order $k$, and
\item[{\bf C.2}] There is a holomorphic function $f$ verifying that
 \begin{equation}\label{ecefe}
 df = \left( c + (a+\ep b) g + \ep \bar{c} g^2\right)\omega 
 \end{equation} 
and that $\omega dg /f$ is holomorphic,
\end{enumerate}    
such that the immersion can be expressed as 
 \begin{equation}\label{recup}
\psi=F\Omega F^*:\Sigma\flecha\L^4. 
 \end{equation}
Here $F:\Sigma\flecha\SL$ is a meromorphic curve satisfying 
 \begin{equation}\label{ecF}
 F^{-1} dF= \left(\begin{array}{cc} 0 & (a +\ep \bar{c} g)\omega\\ dg /f& 
0\\\end{array}\right) =:\cA 
 \end{equation}
 and 
$\Omega:\Sigma\flecha {\rm Herm} (2)$ verifies the differential equation 
\begin{equation}\label{difequation}
{\rm d}\Omega + \mathcal{A} \Omega +\Omega \cA^*= \left(\begin{array}{cc} 
2\varepsilon \Re \left(g \overline{f} \omega\right)& 
\left(1-\varepsilon|g|^2\right) \omega\\ \left(1-\varepsilon|g|^2\right) 
\overline{\omega} & 0\\\end{array}\right). 
\end{equation} 
Conversely, let $\Sigma$ be a simply connected Riemann surface, 
$\varepsilon =\pm 1$, and consider a meromorphic function $g$ and a holomorphic $1$-form 
$\omega$ on $\Sigma$ satisfying {\bf C.1} and {\bf C.2} for some constants $a,b\in \R$, $c\in \C$. 
Then there exist a meromorphic curve $F:\Sigma\flecha\SL$ satisfying 
\eqref{ecF}, and a solution $\Omega:\Sigma\flecha{\rm Herm} (2)$ to the 
system \eqref{difequation}. Moreover, the map 
$\psi:\Sigma\flecha\L^4$ given by  \eqref{recup} is a surface of Bryant 
type in $\L^4$ for which $\varepsilon$ is the sign of its Gauss curvature. 
\end{teo}

\begin{remark}
The only surfaces of Bryant type in $\L^4$ that are flat are those in 
Remark \ref{fl}. So, it is not restrictive to assume in the representation 
theorem that the surfaces of Bryant type are non-flat. 
\end{remark}

\begin{proof}
Let $\p$ be a non-flat simply connected surface of Bryant type. As $Q$ is 
a non-zero holomorphic $2$-form, $\Sigma$ cannot be the Riemann sphere, 
and so we may choose a global holomorphic coordinate $z$ on the Riemann 
surface $\Sigma$. Following the notations of Lemma \ref{flatbundle} we get 
that $q=\tilde{p}-p$ is holomorphic and $f$ as in \eqref{masalfa} is 
meromorphic. On the other hand it is easy to check that 
$d\sigma^2=\varepsilon K ds^2$ is a pseudometric on $\Sigma$ of constant 
curvature $-\varepsilon$. As 
$\Sigma$ is simply connected, by the Frobenius theorem there exists a
meromorphic function $g$ (holomorphic with $|g|<1$ if 
$\varepsilon =1$) on $\Sigma$ such that (see \cite{Bry,GMM2,GaMi})
$$\varepsilon K ds^2 = \frac{4|dg|^2}{\left(1-\varepsilon 
 |g|^2\right)^2}.$$
Now, since from \eqref{defalfa} and \eqref{masalfa} we know that 
$$K=\varepsilon\frac{4|f|^2|q|^2}{\landa^2},$$ it follows that
 \begin{equation}\label{coefdos}
 \frac{|f|^2|q|^2}{\landa}= \frac{|g_z|^2}{\left(1-\varepsilon 
 |g|^2\right)^2}.
 \end{equation}
Thus \begin{equation}\label{coefuno} 
 \landa= \left|\frac{f q}{g_z}\right|^2\left(1-\varepsilon 
 |g|^2\right)^2.
 \end{equation}
In this way $\omega= f Q/dg$ is a meromorphic $1$-form on 
$\Sigma$, and it has no poles. Note that $Q= \omega dg /f$, so this 
quantity defines a holomorphic quadratic differential on $\Sigma$. 
Besides, the zeros of $\omega$ of order $2k$ must trivially coincide with 
the poles of $g$ of order $k$. 

Since from Lemma \ref{conformaleta} we know that the hyperbolic Gauss map 
$G=[\eta+\etat]:\Sigma\flecha\S_{\8}^2$ is conformal, there exist 
holomorphic functions $A,B:\Sigma\flecha\C$ and a positive real function 
$\mu:\Sigma\flecha \R^+$ such that 
 \begin{equation}\label{frn}
\eta+\etat=\mu \left(\begin{array}{cc} A\bar{A} & A\bar{B}\\ \bar{A}B& 
B\bar{B}\\ 
\end{array}\right). 
 \end{equation}
Thus $\esiz (\eta+\etat)_z,(\eta+\etat)_{\bar{z}}\esde=\frac{1}{2} 
|AB_z-BA_z|^2 \mu^2$, and since from \eqref{frame} and \eqref{matrices2} 
it is obtained 
 \begin{equation}\label{recuppsi}
(\eta+\etat)_{\bar{z}} = \frac{2\overline{q}}{\landa} \psi_z, 
 \end{equation}
we get by means of \eqref{coefuno} 
 \begin{equation}\label{liouville}
|A dB-B dA|^2 \mu^2 = \frac{4|dg|^2}{|f|^2\left(1-\varepsilon 
|g|^2\right)^2}. 
 \end{equation}
Besides, from \eqref{liouville} we see that
$$\left| \frac{dg}{f (A dB - BdA)}\right| = \frac{1}{2} \mu (1 -\varepsilon |g|^2 
)>0.$$ So, $dg /(f(AdB -BdA))$ never vanishes, and all its poles are of 
even order. This ensures the existence of a meromorphic function 
(holomorphic if $\varepsilon\neq -1$) $S$ verifying $$ S^2= \frac{dg}{f(A 
dB - B dA)}.$$ If we now choose $C=AS$, $D=BS$ we find that 
$[(A,B)]=[(C,D)]$ and $C dD -D dC = dg /f$. Thus, by substituting 
$(A,B)$ with $(C,D)$ and $\mu$ with $\r$ so that 
$$\eta+\etat=\r \left(\begin{array}{cc} C\bar{C} & C\bar{D}\\ \bar{C}D& 
D\bar{D}\\ \end{array}\right),$$ equation \eqref{liouville} turns into 
 \begin{equation}\label{liouvmod}
\r = \frac{2}{1-\varepsilon |g|^2}. 
 \end{equation}Now consider the meromorphic curve $F:\Sigma\flecha\SL$ 
 \begin{equation}\label{solefej}
F=\left(\begin{array}{cc} C & f dC /dg\\ D& f dD / dg\\ 
\end{array}\right). 
 \end{equation}
Then there exists a meromorphic $1$-form $\vartheta$ on $\Sigma$ such that 
 \begin{equation}\label{efemenosuno}
 F^{-1} dF= 
\left(\begin{array}{cc} 0 & \vartheta\\ d g /f& 0\\\end{array}\right). 
\end{equation}
Moreover, 
 \begin{equation}\label{mincurv}
 \eta+\etat=F \left(\begin{array}{cc}\r & 0\\ 0&0\\ 
\end{array}\right)F^*
 \end{equation}
and from there, \eqref{efemenosuno} and \eqref{liouvmod}, we have
$$(\eta+\etat)_{\bar{z}}=F 
\left(\def\arraystretch{1.5}\begin{array}{cc} 
 \displaystyle\frac{2\varepsilon g\overline{g_z}}{\left(1-\varepsilon |g|^2\right)^2}
 & \displaystyle\frac{2\overline{g_z}}{\overline{f}\left(1-\varepsilon |g|^2\right)}\\ 0&0\\ 
\end{array}\right)F^*.$$ Once here we recall \eqref{recuppsi} and the fact 
that $q$ does not vanish identically to obtain from the above expression 
that 
 \begin{equation*}
\psi_z = \frac{\landa}{2\overline{q}}F 
\left(\def\arraystretch{1.5}\begin{array}{cc} 
 \displaystyle\frac{2\varepsilon g\overline{g_z}}{\left(1-\varepsilon |g|^2\right)^2}
 & \displaystyle\frac{2\overline{g_z}}{\overline{f}\left(1-\varepsilon |g|^2\right)}\\ 0&0\\ 
\end{array}\right)F^*.
 \end{equation*} 
In addition, using \eqref{coefdos} we obtain the final expression for 
$\psi_z$,  
 \begin{equation}\label{psizeta}
\psi_z = F \left(\def\arraystretch{1.5}\begin{array}{cc} 
\displaystyle\frac{\varepsilon |f|^2 q g}{g_z}& 
\displaystyle\frac{\left(1-\varepsilon|g|^2\right) q f }{g_z} \\ 0& 
0\\\end{array}\right) F^*. 
 \end{equation} 
Finally, let us note that as $\SL$ acts on $\L^4$ as the connected 
component of the identity in its isometry group, the immersion 
$\p$ can be expressed as 
$\psi=F\Omega F^*:\Sigma\flecha\L^4$ for an adequately chosen intermediate matrix 
$\Omega:\Sigma\flecha {\rm Herm }(2)$. It comes 
clear from \eqref{psizeta} that $\Omega$ is a solution of the differential 
system \eqref{difequation}. Next, note that by differentiation of 
\eqref{psizeta} with respect to $\bar{z}$ and noting that 
$\psi_{z\bar{z}}$ is real, we see that the data $f,g,\omega, \vart$ must verify 
$$\ep g \omega \left(\overline{df} - \bar{g} \bar{\vart}\right) + \omega \bar{\vart} = 
\ep \bar{g} \bar{\omega }\left(df - g \vart\right) + \bar{\omega} \vart.$$ 
Thus, by Lemma \ref{apen} in the Appendix, and as $g$ is not constant 
(otherwise $\psi$ would be flat), we obtain the existence of constants 
$a,b\in \R$ and $c\in \C$ such that
 \begin{equation}\label{latita}
 \vart = (a +\ep \bar{c} g)\omega, \hspace{1cm} df =  \left( c + 
(a+\ep b) g + \ep \bar{c} g^2\right) \omega. 
 \end{equation}
Particularly, $f,\vart$ are holomorphic, and {\bf C.2} and \eqref{ecF} 
hold. Thus the proof of the first part of the theorem is complete. 

For the converse, we start with the Weierstrass data $(g,\omega)$ and 
constants $a,b\in \R$, $c\in \C$ verifying 
${\bf C.1}$ and ${\bf C.2}$ on a simply connected Riemann surface 
$\Sigma$. Then the system \eqref{ecF}  has a (possibly multivalued) 
solution $F:\Sigma\flecha \SL$ of the form \eqref{solefej}, where $C,D$ 
are linearly independent solutions of (see \cite{GMM0}) 
 \begin{equation}\label{singuregu}
 Z'' - \frac{(g'/f)'}{g'/f} Z' - (a +\ep \bar{c} g) q Z =0, \hspace{0.5cm} 
\left(' = \frac{d}{dz}\right),
 \end{equation}
being $Q= q(z) dz^2 = \omega dg /f.$

As the meromorphic $1$-form in \eqref{ecF} has its poles at the poles of 
$g$, we see that the solution $F$ is locally well defined and holomorphic away from the 
poles of $g$.

Let now $z_0 \in\Sigma$ be a pole of $g$ of order $k\geq 1$, and let 
$\delta \geq 0$ denote the order of the zero of $f$ at $z_0$ (possibly $\delta =0$). From {\bf 
C.2} it is clearly seen that if $c\neq 0$, then $\delta \in \{0,1\} $, 
while if $c=0$ then $\delta =0$ (otherwise we would have $\delta = k+1$ by 
{\bf C.2}, which would contradict that $\omega dg /f$ has no poles).

It is straightforward that $(g'/f)' /(g'/f)$ has a simple pole at $z_0$, 
of residue $-(k+\delta +1)$. Besides, if $h= (a +\ep \bar{c} g)q$, it 
follows directly that $h$ is holomorphic at $z_0$ if $c=0$, it has a 
simple pole at $z_0$ if $c\neq 0$ and $f(z_0)\neq 0$, and has a pole of 
order two at $z_0$ otherwise. From this, a simple calculation shows that 
$$h_{-2} = \lim_{z\to z_0} (z-z_0)^2 h(z)= \delta k.$$ Hence, the 
differential equation \eqref{singuregu} has a regular singularity at $z_0$ 
(see \cite{CoLe}). Its indicial equation is $$\landa^2 + (k +\delta)\landa 
+ k \delta =0,$$ that has the integer roots $-k$ and $-\delta$. Therefore, 
both $C,D$ are single valued meromorphic functions on $\Sigma$, and the 
orders of their poles at $z_0$ are $k$ and $\delta$ (see \cite{GMM0}). 
Particularly, any pole of $C$ or $D$ of order $l$ is located at a poles of 
$g$ of order $\geq l$.

Once here, we have ensured the existence of a meromorphic solution 
$F:\Sigma\flecha \SL$ to \eqref{ecF} of the form \eqref{solefej}. If we now set $\varrho:\Sigma\flecha 
[0,+\8)$ as \eqref{liouvmod}, the map 
 \begin{equation}\label{mincur2}
N= F\left(\begin{array}{cc}\r & 0\\ 0&0\\ 
\end{array}\right)F^* 
 \end{equation}
has a finite value at every point, due to the previous analysis regarding 
the poles of $C,D$. Thus we have a map 
$N:\Sigma\flecha \N^3$. The same argument shows that the $1$-forms $$\phi 
dz = F \left(\def\arraystretch{1.5}\begin{array}{cc} \varepsilon g \bar{f} 
\omega & \left(1-\varepsilon|g|^2\right) \omega \\ 0& 
0\\\end{array}\right) F^*, \hspace{0.3cm} \widetilde{\phi} d\bar{z} = F 
\left(\def\arraystretch{1.5}\begin{array}{cc} \varepsilon \bar{g} f 
\bar{\omega} & 0 \\ \left(1-\varepsilon|g|^2\right) \bar{\omega} &
0\\\end{array}\right) F^*$$ take finite values at all points. But now, 
noting that $\phi^* = \widetilde{\phi}$ and that $\phi_{\bar{z}} = 
\widetilde{\phi}_z$, we can conclude the existence of a map 
$\psi:\Sigma\flecha \He$ such that $\psi_z =\phi$ and $\psi_{\bar{z}} =\widetilde{\phi}$. 
Finally, let us define $\Omega :\Sigma \flecha \He$ as $\Omega = F^{-1} 
\psi (F^{-1})^*$, whose entries may take infinite values at some points. 
Then $\psi =F\Omega F^*$ and $\Omega$ is trivially a solution to the 
differential system \eqref{difequation}.

At last, as by differentiation of \eqref{psizeta} we have that 
$\psi_{z\bar{z}}$ is collinear with $N$, we conclude that $\psi$ is a marginally trapped 
surface. Now, deriving the right hand side of \eqref{mincur2} with respect 
to $z$ we get that $N$ is parallel, that is $\psi$ has flat normal bundle. 
And as $\esiz d\psi,d\psi\esde = (1-\ep |g|^2 )^2 |\omega|^2$, $\psi$ is 
regular and of Bryant type, and the proof is finished. 
\end{proof} 

\begin{remark}
The mean curvature vector of the Bryant-type surface $\psi:\Sigma\flecha 
\L^4$ verifies the relation $2 \psi_{z\bar{z}} = \landa {\bf H} = 2 E 
(\eta + \etat).$ This indicates by differentiation of \eqref{psizeta} that 
 \begin{equation}\label{eland}
 \frac{2 E}{\landa} = \frac{a+b|g|^2 +2\ep \Re (\bar{c} g)}{1-\ep |g|^2}.
 \end{equation}
Thus an explicit expression for 
${\bf H}$ in terms of the Weierstrass data is obtained from \eqref{mincu}, 
\eqref{eland} and \eqref{mincurv}.
\end{remark}

\begin{remark}
The hyperbolic Gauss map $G:\Sigma\flecha \C\cup \{\8\} $ of a surface of 
Bryant type is a geometric concept, and thus is uniquely determined at 
every point. However, this is not the case for the other basic meromorphic 
data of a Bryant surface. First, note that as $Q$ depends on the chosen 
frame of the normal bundle, it is defined up to the change $Q \to 
e^{\alfa} Q$, where 
$\alfa \in \R$ is the constant hyperbolic angle relating two different frames. 
Therefore, the function $f$ is by definition defined up to 
$f\to e^{-\alfa + i \beta} f$, $\beta\in \R$. Besides, the meromorphic 
function $g$ is unique up to isometries of the $2$-sphere 
$\S^2 \equiv \C\cup \{\8\} $ if $\ep =-1$, and up to isometries of the 
Poincaré disk $\H^2 \equiv \D$ if $\ep =1$. Thus, $g$ is unique up the 
change $$g \to \frac{\tau g +\ep \bar{\gamma}}{\gamma g + \bar{\tau}}, 
\hspace{0.5cm} |\tau|^2 - \ep |\gamma|^2 =1.$$ Noting now that 
$Q=\omega dg / f$, the above comments show that $\omega$ is defined up to 
$$\omega \to e^{i\beta} (\gamma g + \bar{\tau} )^2 \omega.$$ 
\end{remark}

To close this section, we shall relate the hyperbolic Gauss map 
$G:\Sigma\flecha \C\cup \{\8\} $ to the Weierstrass data of a Bryant-type 
surface in $\L^4$. First, observe that in \eqref{frn} we may choose 
$A=1$, $B=G$. Then, with these choices we end up with the 
formula $\psi = F\Omega F^*$ where 
 \begin{equation*}
F= \left( \begin{array}{cc} C & f dC / dg \\ D & f dD / dg 
\end{array} \right) ,\hspace{0.3cm} C =\sqrt{dg / (fdG)}, \ D = G\sqrt{g / 
(fdG)}. 
 \end{equation*}
This formula extends a result by Small \cite{Sma} for mean curvature one 
surfaces in $\H^3$ (see also \cite{GMM2}). 
 
An alternative relation between $G$ and the Weierstrass data relies in the 
concept of \emph{Schwarzian derivative} 
$\{h,z\} $ of a meromorphic function $h$: $$\{h,z\} = 
\left(\frac{h''}{h'}\right)' - \frac{1}{2}\left(\frac{h''}{h'}\right)^2, 
\hspace{1cm} \left( '=\frac{d}{dz}\right).$$ We get then by \cite[Eq. 
(28)]{GMM0} and \eqref{efemenosuno} that the following relation holds on 
any surface of Bryant type in $\L^4$: 
 \begin{equation*}
\{G,z\} dz^2= \left(\left(\frac{(g'/f)'}{g'/f}\right)' -\frac{1}{2} 
\left(\frac{(g'/f)'}{g'/f}\right)^2 \right) dz^2 - (a +\ep \bar{c} g) Q. 
 \end{equation*}
We remark that this formula extends an important equation due to Umehara 
and Yamada \cite{UmYa2} in the context of mean curvature one surfaces in 
$\H^3$.

\section{Examples}

{\bf Representation of CMC-$r$ surfaces:} next, we show that the conformal 
representation in Theorem \ref{rep} generalizes the Bryant representation 
\cite{Bry} for surfaces with $H=r$ in $\H^3 (-r^2)$, as well as the 
Aiyama-Akutagawa one \cite{AiAk} for spacelike surfaces with 
$H=r$ in 
$\S_1^3 (r^2)$. For this, we shall use the unified notation $\M^3 (\ep 
r^2)$ to denote $\H^3 (-r^2)$ for 
$\ep =-1$ and $ \S_1^3 (r^2)$ for $\ep =1$.

Let $\psi:\Sigma\flecha \M^3(\ep r^2)\subset \L^4$ be a simply connected 
(spacelike) CMC-$r$ surface, and let $\eta_1$ be its unit normal in 
$\M^3(\ep r^2)$. Then $\{\eta_1,\eta_2 :=r \psi\}$ is a parallel 
orthonormal frame in the normal bundle of $\psi$ in $\L^4$. So, using the 
notations of the first two sections, it follows directly that $\ep Q = 
\esiz \psi_{zz},\eta_1 \esde \, dz^2$. Therefore $|f|=1$ on $\Sigma$ and 
as $f$ is defined up to constant rotations (note that in this case we are 
working with a uniquely determined frame in the normal bundle), we may 
assume that $f=1$.

With all of this, the differential system \eqref{difequation} can be 
explicitly solved under the condition $\det (\Omega )= -\ep /r^2$, to 
obtain $$\Omega = \frac{1}{r} \left(\def\arraystretch{1.1} 
\begin{array}{cc} -\ep & \ep \bar{g} \\ \ep g & 1 - \ep |g|^2 \end{array} \right) = 
\frac{1}{r} \left(\def\arraystretch{1.1} 
\begin{array}{cc} 0 & i \\ i & -i g \end{array} \right)
\left(\def\arraystretch{1.1} 
\begin{array}{cc} 1 & 0 \\ 0 & -\ep \end{array} \right)
 \left(\def\arraystretch{1.1} \begin{array}{cc} 0 & -i \\ 
-i & i \bar{g} \end{array} \right) .$$ Finally, we derive with respect to 
$z$ the identity $\psi = F\Omega F^*$ and compare it with \eqref{psizeta} 
to deduce that $\vartheta = r \omega$ for the holomorphic $1$-form 
$\vart$ in \eqref{latita}. In conclusion, the Bryant surface 
$\psi:\Sigma\flecha \M^3 (\ep r^2)$ is recovered as $$\psi = \frac{1}{r} \mathcal{B} 
\left(\def\arraystretch{1.1} 
\begin{array}{cc} 1 & 0 \\ 0 & -\ep \end{array} \right) \mathcal{B}^*, 
\hspace{0.4cm} \mathcal{B}= F \left(\def\arraystretch{1.1} 
\begin{array}{cc} 0 & i \\ i & -i g \end{array} \right).$$ Here 
$\mathcal{B}:\Sigma\flecha \SL$ is a null holomorphic curve (i.e. $\det 
(d\mathcal{B})=0$), and $F:\Sigma\flecha \SL$ verifies \eqref{ecF} for 
$a=r$, $c=0$. Thus, the Bryant representation in \cite{Bry} and the Bryant-type representation
in \cite{AiAk} are recovered.

{\bf The Weierstrass representation:} now, we prove that the classical 
Weierstrass representation for minimal surfaces in $\R^3$ and its analogue 
for maximal surfaces in $\L^3$ are also included in Theorem \ref{rep}. In 
order to do so, we fix the notation 
$\R_{\ep}^3 $ to denote $\R^3 \equiv x_0 =0 \subset \L^4$ if $\ep =-1$, and 
$\L^3 \equiv x_3 =0 \subset \L^4$ if $\ep =1$.

Let $\psi:\Sigma\flecha \R_{\ep}^3 \subset \L^4$ be a minimal (or maximal) 
surface in $\R_{\ep}^3$, with unit normal $\eta_1$. Then $\{\eta_1, 
\eta_2=\frac{1}{2}(1-\ep,0,0,1+\ep)\}$ is a parallel orthonormal frame of 
the normal bundle of 
$\psi$ in $\L^4$, and arguing as above we get that $f=1$. Moreover, as ${\bf H}=0$, by \eqref{eland}
we have $a=b=c=0$. Thus the differential equation \eqref{ecF} can be 
explicitly integrated, and we obtain a solution as 
\begin{equation}\label{Fminimal} 
F=\left( 
\def\arraystretch{1.1} 
\begin{array}{cc} 1 & 0 \\ g & 1 \end{array} \right) :\Sigma \flecha \SL .
\end{equation} 
Let $\Omega:\Sigma\flecha \He$ be the solution to \eqref{difequation}, 
which is in this case as 
 \begin{equation}\label{raao}
\Omega_z + \left( 
\def\arraystretch{1.1} 
\begin{array}{cc} 0 & 0 \\ g_z & 0 \end{array} \right) \Omega = \left( 
\def\arraystretch{1.1} 
\begin{array}{cc} \ep g \omega & (1 -\ep |g|^2 ) \omega  \\ 0 & 0 \end{array} 
\right).
 \end{equation}
If we write $$\Omega = \left( 
\def\arraystretch{1.1} 
\begin{array}{cc} U & V \\ \bar{V} & W \end{array} \right) ,$$ then we see 
from \eqref{raao} that $$U = 2\ep \Re \int g \omega, \hspace{0.4cm} (U g + 
\bar{V})_z = \ep g^2 \omega, \hspace{0.4cm} (U \bar{g} + V)_{\bar{z}} = 
\omega.$$ Now, as $\psi (\Sigma)\subset \R_{\ep}^3$ it holds $$\psi = 
F\Omega F^* = \left( 
\def\arraystretch{1.1} 
\begin{array}{cc} U & U \bar{g} +V \\ U g + \bar{V} & \ep U \end{array} \right) 
.$$ So, by the above computations we finally obtain the Weierstrass 
representation: $$\psi = \Re \int \left( (1+ \ep ) g, 1 +\ep g^2, -i(1- 
\ep g^2), (1- \ep) g \right) \omega.$$ 

{\bf Analytic deformation of surfaces:} we shall show now that, with the 
above notations, CMC-$r$ surfaces in $\M^3 (\ep r^2)$ can be analytically 
and isometrically deformed to minimal (or maximal) surfaces in 
$\R_{\ep}^3$ as $r$ tends to zero. This result was obtained in the case $\ep =-1$ by Umehara 
and Yamada \cite{UY}.

To do so, we consider first the \emph{translated spaces} $$\widetilde{{\rm 
M}}^3 (\ep r^2)= \left\{p - \frac{1}{2r}(1-\ep, 0,0,-1-\ep) : p \in \M^3 
(\ep r^2)\right\} ,$$ where we will assume that $\M^3(r^2) = \{x\in \S_1^3 
: x_3<0 \} $ is a de Sitter half-space. Our deformation process will rely 
on the fact that, as $r\to 0$, the spaces $ \widetilde{{\rm M}}^3 (\ep 
r^2) $ converge to $\R_{\ep}^3$. 

Let $\psi_{r_0} :\Sigma\flecha \M^3 (\ep r_0^2)$ be a simply connected 
CMC-$r_0$ surface, and choose $z_0\in \Sigma$. Let now $(g,\omega)$ be its 
Weierstrass data, and suppose without loss of generality that $g(z_0)=0$. 
Then, for every $r>0$ there exists a unique (up to rigid motions) CMC-$r$ 
surface in $\M(\ep r^2)$ $\psi_r :\Sigma\flecha \M(\ep r^2)$ that has the 
same Weierstrass data $(g,\omega)$, and so it is isometric to the original 
immersion $\psi_{r_0}$. Now, we know that
$\psi_r = \frac{1}{r} F_r \Delta F_r^* 
$, where $F_r\in \SL$, 
$$F_r^{-1} dF_r = \left(\begin{array}{cc} 0 & r\omega \\ dg & 0 
\end{array} \right), \hspace{0.6cm} F_r (z_0)={\rm Id}, \hspace{0.6cm} 
\Delta = \left(\def\arraystretch{1.1} 
\begin{array}{cc} -\ep & \ep \bar{g} \\ \ep g & 1 - \ep |g|^2 \end{array} 
\right).$$ Consider now the translated immersions $X_r :\Sigma\flecha 
\widetilde{{\rm M}}^3 (\ep r^2)$ given by $$X_r= \psi_r - \frac{1}{2r} 
(1-\ep, 0,0,-1-\ep) = \frac{1}{r} \left( F_r \Delta F_r^* - 
\left(\def\arraystretch{1.1} 
\begin{array}{cc} -\ep & 0 \\ 0 & 1 \end{array} 
\right) \right).$$ As the family $\mathcal{A}_r := F_r^{-1} dF_r $ is real 
analytic with respect to $r\in \R$, the family $F_r:\Sigma\flecha \SL$ is 
also real analytic with respect to $r$. So, using that $F_0$ is given by 
\eqref{Fminimal}, it is easy to see that 
$$F_r \Delta F_r^* - \left(\def\arraystretch{1.1} 
\begin{array}{cc} -\ep & 0 \\ 0 & 1 \end{array} 
\right) = a_1 (z,\bar{z}) r + o (r).$$ This assures that the family of 
surfaces $X_r :\Sigma\flecha \L^4$ is real analytic with respect to 
$r\in \R$. Furthermore, $X_0:\Sigma \flecha \R_{\ep}^3$ is given by $$X_0 = 
a_1(z,\bar{z})= \left. \frac{\parc}{\parc r} \right|_{r=0} F_r \Delta 
F_r^*.$$ Finally, $X_0$ has zero mean curvature in $\R_{\ep}^3$. This 
happens because, by $$\frac{\parc^2}{\parc z\parc \bar{z}} (F_r \Delta 
F_r^*) = F_r \left(\def\arraystretch{1.1} 
\begin{array}{cc} r^2 |\omega|^2 (1-\ep |g|^2) & 0 \\ 0 & 0 \end{array} 
\right) F_r^* $$ and the analyticity of the family in $r\in \R$, we have 
$$ \frac{\parc^2 X_0}{\parc z\parc \bar{z}} = 
\frac{\parc^2}{\parc z\parc \bar{z}} \left(\left.\frac{\parc}{\parc r} 
\right|_{r=0} F_r \Delta F_r^* \right) = \left.\frac{\parc}{\parc r} 
\right|_{r=0} \left(\frac{\parc^2}{\parc z\parc \bar{z}} (F_r \Delta 
F_r^*) \right) =0.$$ 

Therefore, we conclude that the CMC-$r$ surfaces of $\H^3(-r^2)$ (resp. 
$\S_1^3 (r^2)$) can be perturbed in an analytic and isometric way to minimal surfaces in 
$\R^3$ (resp. to maximal surfaces in $\L^3$) as $r$ approaches to zero.
 
\vspace{0.1cm}
 
{\bf New families of complete examples:} there exist many complete 
surfaces of Bryant type in $\L^4$ with non-parallel mean curvature, that 
can be constructed by means of the representation formula. 

To see this, let us choose $c=0$ in Theorem \ref{rep}. Then $df= (a+\ep b) 
g \omega$. Therefore, the condition {\bf C.2} in this case just asks for 
the existence of a nowhere-zero primitive $\int g \omega$ of the 
holomorphic 
$1$-form $g\omega$ on $\Sigma$, a condition that always holds locally. 
Moreover, if $a+ \ep b \neq 0$ the example has non-parallel mean 
curvature, and if $a=-\ep b\neq 0$ the surface lies in some affine 
hyperbolic $3$-space or de Sitter $3$-space in $\L^4$ and is a Bryant 
surface there (see Remark \ref{rem7}). Finally, if $a=-\ep b=0$, the 
surface has zero mean curvature in some Euclidean or Lorentzian affine 
$3$-space of $\L^4$, again by Remark \ref{rem7}.

Many of the examples with $c=0$ and non-parallel mean curvature are 
complete. For instance, if $(g,\omega)$ are the Weierstrass data of a 
complete minimal surface in $\R^3$ lying in a halfspace with horizontal 
boundary, then its third coordinate $\Re \int g\omega$ is non-surjective, 
and thus the above condition {\bf C.2} holds. So, $g,\omega$ together with 
$\ep =-1$ and 
$a,b\in \R$, $a+\ep b \neq 0$, generate a complete surface of Bryant type 
in $\L^4$. 

As a closing remark, we indicate that if $a=c=0$ and $b\neq 0$, the 
resulting class of surfaces admit an integral representation quite 
analogous to the Weierstrass representation of minimal and maximal 
surfaces. Indeed, in that situation we have $\vart=0$ in \eqref{latita}, 
so a solution $F:\Sigma\flecha \SL$ to \eqref{ecF} is obtained by 
substituting $g$ in \eqref{Fminimal} by 
$\int dg /f$. With this, the system \eqref{difequation} can be integrated 
much in the same way that we did for minimal and maximal surfaces in 
\eqref{raao}. We do not write the final expressions explicitly, as they 
are straightforward but rather lengthy.

\section{Completeness}

The well known Calabi-Bernstein theorem \cite{Cal} asserts that the only 
complete maximal surfaces in $\L^3$ are spacelike planes. Analogously, 
every complete spacelike CMC-$r$ surface in $\S_1^3 (r^2)$ must be a flat 
totally umbilic example, obtained as the intersection of $\S_1^3 (r^2)$ 
with a degenerate vector hyperplane of $\L^4$ \cite{Aku,Ram}. We remark 
that both maximal surfaces in $\L^3$ and CMC-$r$ surfaces in $\S_1^3 
(r^2)$ have non-negative curvature, $K\geq 0$. We start this section with 
a simultaneous generalization of these two Bernstein-type theorems: 
 \begin{cor}\label{12}
Every complete surface of Bryant type in $\L^4$ with non-negative 
curvature is a flat surface lying in a degenerate hyperplane of $\L^4$, as 
described in Remark \ref{fl} 
 \end{cor} 
\begin{proof}
Given a surface of Bryant type $\psi:\Sigma\flecha \L^4$ with non-negative 
curvature, we obtain that $\varepsilon =1$, and so $\esiz d\psi,d\psi 
\esde = (1 - |g|^2)^2 |\omega|^2 \leq |\omega|^2$. Thus $|\omega|^2$ is a 
flat metric, which is complete because so is $\psi$. Therefore, the 
Riemann surface $\Sigma$ must be parabolic, and as $|g|<1$ we obtain that 
$g$ is constant. Therefore the Hopf differential $Q$ vanishes identically on $\Sigma$, and the 
surface must be flat and lie in a degenerate hyperplane of $\L^4$, by 
Remark \ref{fl}. 
\end{proof}

In minimal surface theory, as well as in Bryant surface theory, the study 
of the complete examples of finite total curvature has been widely 
developed. Here, we recall that a surface $\Sigma$ with non-positive 
curvature $K \leq 0$ has finite total curvature provided $$\int_{\Sigma} K 
dA > -\8,$$ where $dA$ is the area element of the surface. 

Our next result shows that, even though there are many complete simply 
connected surfaces of Bryant type with non-parallel mean curvature, none 
of them has finite total curvature.
 \begin{teo}
Let $\psi:\Sigma\flecha \L^4$ be a non-flat complete simply-connected 
surface of Bryant type with finite total curvature. Then $\psi(\Sigma)$ 
lies in some Euclidean or hyperbolic $3$-space of $\L^4$, and its 
Weierstrass data are given by $$ g(z)= \frac{P_1(z)}{P_2(z)}, \hspace{1cm} 
\omega = P_2(z)^2 dz , \hspace{0.5cm} c = a+\ep b =0,$$ where 
$P_1(z),P_2(z):\C\flecha \C$ are polynomials with no common zeros. 
 \end{teo}
\begin{proof}
Since $\psi$ is a non-flat immersion then $\ep=-1$, from Corollary 
\ref{12}. On the other hand, using that $\psi$ has finite total curvature, 
$\Sigma$ must be parabolic, that is, we can assume $\Sigma=\C$, and its 
Weierstrass data $(g,\omega)$ are meromorphic on $\C\cup\{\infty\}$ (see 
\cite{Oss}). 

As $\omega$ is holomorphic on $\Sigma=\C$, there exist 
$Q_1(z),Q_2(z),Q_3(z):\C\flecha \C$ polynomials, $Q_1(z),Q_2(z)$ without
common factors, such that $$ g(z)= \frac{Q_1(z)}{Q_2(z)}, \hspace{1cm} 
\omega = Q_3(z) dz. $$ Observe that from {\bf C.1}, the zeros of $g$ of 
order $k$ correspond to the zeros of order $2k$ of the holomorphic 
$1$-form $g^2\omega$. So, 
$g^2\omega=A Q_1(z)^2\,dz$ for a non-zero complex constant $A$. Then, 
$Q_3(z)=A Q_2(z)^2$ and we can write $ g(z)=P_1(z)/P_2(z),\, \omega = 
P_2(z)^2 dz$ for 
$P_1(z)=\sqrt{A}Q_1(z),\,P_2(z)=\sqrt{A}Q_2(z)$. 

Let $R_1, R_2$ be the degrees of $P_1(z),P_2(z)$, respectively. Since $g$ 
is unique up to isometries of $\C\cup\{\infty\}$, we can suppose that 
$g(\infty)=\infty$, that is, $R_1>R_2$. Consequently 
\[
{\rm degree}(\omega)<{\rm degree}(g\omega)<{\rm degree}(g^2\omega). 
\]
Moreover, from (\ref{ecefe}), $f(z)$ is a polynomial such that ${\rm 
degree}(df)={\rm degree}(g^2\omega)$ if $c\neq 0$ or ${\rm 
degree}(df)={\rm degree}(g\omega)$ if $c=0$ and $a- b\neq 0$. Hence, if 
$c\neq 0$ or $a-b\neq 0$ it follows that ${\rm degree}(df)\geq {\rm 
degree}(g\omega)=R_1+R_2$, and so 
$${\rm degree}(f)\geq R_1+R_2+1>R_1+R_2-1\geq {\rm degree}(\omega dg).$$
But the last inequality implies that $\omega dg/f$ is not holomorphic, 
which is a contradiction. Thus we can conclude that $c=0$ and $a- b= 0$, 
which means by Remark \ref{rem7} that $\bf{H}$ is parallel, and $\psi$ 
lies in some Euclidean or hyperbolic $3$-space of $\L^4$. This concludes 
the proof. 
\end{proof}

\section*{Appendix}

In searching an adequate family which generalizes simultaneously the 
theories of minimal surfaces in $\R^3$ and Bryant surfaces in $\H^3$, a 
natural hypothesis is to ask the mean curvature vector to be parallel in 
the normal bundle, that is, $\nabla^{\perp} {\bf H} \equiv 0$. However, 
this condition is too strong, as it does not generate any new example: 
 \begin{pro}
Any spacelike surface in $\L^4$ with parallel mean curvature which is not
maximal in $\L^4$ must lie in a totally umbilical affine hyperquadric 
$\Q^3$ of $\L^4$.
 \end{pro} 
\begin{proof}
Let $\p$ be a non maximal spacelike surface with parallel mean curvature 
${\bf H}$. Then it holds that $\langle{\bf H},{\bf H}\rangle=A$ with 
$A\in\R$. 

If $A=0$, then $\psi$ is a marginally trapped surface, and we can choose 
$\nu\in T^\bot(\Sigma)$ such that $\langle\nu,\nu\rangle=0$ and 
$\langle\nu,{\bf H}\rangle=C>0$. Then it is easy to check that 
$\{\eta=(\nu+{\bf H})/(2C),\widetilde{\eta}=(\nu-{\bf H})/(2C)\}$ is an
orthonormal frame of $T^\bot(\Sigma)$, and that both 
$\eta,\widetilde{\eta}$ are parallel. This implies that the normal 
curvature vanishes identically, 
$R^\bot\equiv 0$. The result follows then by Remark \ref{rem7}.

If $A\neq 0$, let us suppose that $A>0$ and so ${\bf H}$ is spacelike (the 
case $A<0$ is analogous). Then we can take $\widetilde{\eta}$ a normal 
timelike vector field such that $\{\eta={\bf H}/A,\widetilde{\eta}\}$ is 
an orthonormal frame of $T^\bot(\Sigma)$. As both 
$\eta,\widetilde{\eta}$ are parallel, $R^\bot\equiv 0$. Since the mean 
curvature vector ${\bf H}$, given by \eqref{laac}, is spacelike, we get 
$\tilde{E}=0$ and $E/\lambda=0$. Now, by the Codazzi and Ricci equations 
in (\ref{GCR}) it follows that $p,\widetilde{p}$ are holomorphic and 
$\widetilde{p}=k p$, $k\in \R$. The proof of this case finishes following 
the argument of Remark \ref{rem7}.

\end{proof}

The following elementary fact is used in the proof of the representation 
theorem. 

\begin{lem}\label{apen}
Let $f_1,f_2,f_3,f_4:\Sigma\flecha \C\cup \{\8\} $ be meromorphic 
functions on a Riemann surface, such that $f_1$ and $f_3$ are linearly 
independent and \begin{equation}\label{realiti} f_1 \bar{f_2} + f_3 
\bar{f_4} = \bar{f_1} f_2 + \bar{f_3} f_4. 
\end{equation} 
Then there exist constants $a,b\in \R$ and $c\in \C$ such that $$f_2 = a 
f_1 + c f_3, \hspace{1cm} f_4 = \bar{c} f_1 + b f_3.$$ 
\end{lem}
\begin{proof}
Given $z_0 \in \Sigma$ a point where the functions $f_i$ have no poles and 
$f_3(z_0) \neq 0$, by differentiation of \eqref{realiti} we get 
 \begin{equation}\label{ap1}
 f_1^{(n)} \bar{f_2} + f_3^{(n)} \bar{f_4} = \bar{f_1} f_2^{(n)} +  \bar{f_3} f_4^{(n)}
 \end{equation}
at $z_0$, for every $n\in \N$. Thus, all derivatives of $f_4$ at $z_0$ are 
a linear combination of the derivatives of $f_1,f_2,f_3$ at $z_0$. This 
shows the existence of $\landa, \mu, \delta \in \C$ with 
 \begin{equation}\label{ap2}
 f_4 = \landa f_1 + \mu f_2 + \delta f_3.
 \end{equation} 
Deriving now \eqref{ap1} with respect to $\bar{z}$ at $z_0$ we obtain 
 \begin{equation}\label{ap3}
 f_1^{(n)} \overline{f_2^{(k)}} + f_3^{(n)} \overline{f_4^{(k)}} = \overline{f_1^{(k)}} f_2^{(n)} +  
 \overline{f_3^{(k)}} f_4^{(n)}
 \end{equation}
for all $k\in \N$. As $f_1,f_3$ are linearly independent, there is some 
$n_0\in \N$ such that $f_1^{(n_0)} + \bar{\mu}  f_3^{(n_0)} \neq 0$ at 
$z_0$. Putting this fact together with \eqref{ap2} and \eqref{ap3} we 
conclude as before that $f_2$ is a linear combination of $f_1$ and $f_3$. 
Therefore, there exist complex constants $a,b,c,e\in \C$ such that 
 \begin{equation}\label{ap4}
 f_2 = a f_1 + c f_3, \hspace{1cm} f_4 = e f_2 + b f_4.
 \end{equation} 
Now, by \eqref{realiti} and \eqref{ap4} it holds $$ (a- \bar{a}) |f_1|^2 + 
(c-\bar{e}) \bar{f_1} f_3 + (e - \bar{c}) f_1 \bar{f_3} + (b-\bar{b}) 
|f_3|^2 =0.$$ By the linear independence of $f_1$ and $f_3$ this indicates 
that $a,b\in \R$ and $e =\bar{c}$, what completes the proof. 
\end{proof}

 \end{document}